\documentclass[12pt]{article}
\usepackage{epic,eepic,epsf,epsfig}
\usepackage{amsmath,enumerate}
\usepackage{amssymb}
\usepackage{amsbsy}

\usepackage{amsfonts}
\newtheorem{theorem}{Theorem}[section]%
\newtheorem{lemma}[theorem]{Lemma}%
\newtheorem{cor}[theorem]{Corollary}%
\newtheorem{question}[theorem]{Question}%

 \def\Omega{\Omega}

\def\f{\noindent}

\newcommand{\qed}{\mbox{\raisebox{0.7ex}{\fbox{}}} \vspace{4truemm}}
\def\demo{\f {\bf Proof.}\hskip10pt}

\begin{document}

\baselineskip 16pt

\title{ \vspace{-1.2cm}
On a finite group in which every maximal subgroup of order divisible by $p$ containing the normalizer of some Sylow subgroup is nilpotent
\thanks{\scriptsize This work was supported in part by NSFC (12571024) and Shandong Provincial Natural Science Foundation, China (ZR2017MA022).
\newline
 \hspace*{0.5cm} \scriptsize $^{\ast\ast}$Corresponding
  author.\newline
 \hspace*{0.5cm} \scriptsize{E-mail addresses:} jtshi2005@pku.org.cn\,(J. Shi),\,\,lifh2026@s.ytu.edu.cn\,(F. Li).}}

\author{Jiangtao Shi\,$^{\ast\ast}$,\,\,Fanghua Li\\
\\
{\small{\em School of Mathematics and Information Sciences,
Yantai University, Yantai 264005, P.R. China}}}

\date{ }

\maketitle \vspace{-.8cm}

\begin{abstract}
Let $G$ be a finite group and $p$ any fixed prime divisor of $|G|$. We prove that if
every maximal subgroup of $G$ of order divisible by $p$ containing the normalizer of some Sylow subgroup is nilpotent then $G$ is a solvable group having
a Sylow tower and $G$ is $q$-nilpotent or $q$-closed for each prime divisor $q$ of $|G|$.

\medskip \f {\bf Keywords:} maximal subgroup; normalizer; Sylow subgroup; nilpotent; $q$-nilpotent; $q$-closed
\end{abstract}

\section{Introduction}\label{s1}

In this paper all groups are assumed to be finite. In {\rm\cite{redei1956}} R$\rm\acute{e}$dei provided a complete characterization of Schmidt groups in which every maximal subgroup is nilpotent. Deng, Meng and Lu {\rm\cite[Theorem 3.1]{meng}} generalized {\rm\cite{redei1956}} and they gave a description of a group in which every maximal subgroup of even order is nilpotent. Let $p$ be any fixed prime divisor of the order of a group $G$.
Shi, Li and Shen {\rm\cite[Theorem 1.1]{shili}} indicated that a non-solvable group $G$ always has a non-nilpotent maximal subgroup of order divisible by $p$. Moreover, Shi and Tian {\rm\cite[Theorem 1.1]{shitian}} obtained a complete characterization of a group $G$ in which every maximal subgroup of order divisible by $p$ is nilpotent.

Considering the coprime action of groups, let $A$ and $G$ be groups such that $A$ acts coprimely on $G$ by automorphisms. {\rm\cite{jiangshao}} and {\rm\cite{shishanxu}} independently characterized the structure of a group in which every maximal $A$-invariant subgroup of order divisible by $p$ is nilpotent from different perspectives using different proof methods. In particular, Shi and Xu {\rm\cite[Theorem 1.3]{shixu2025}}
provided a complete characterization of a group $G$ in which every maximal $A$-invariant subgroup containing the normalizer of some
$A$-invariant Sylow subgroup is nilpotent.

In this paper, as a further generalization of above research, we have the following result, the proof of which is given in Section~\ref{s3}.

\begin{theorem}\ \ \label{th1} Let $G$ be a group and $p$ any fixed prime divisor of $|G|$. If
every maximal subgroup of $G$ of order divisible by $p$ containing the normalizer of some Sylow subgroup is nilpotent. Then

$(1)$ $G$ is solvable;

$(2)$ $G$ has a Sylow tower;

$(3)$ $G$ is $q$-nilpotent or $q$-closed for each prime divisor $q$ of $|G|$.
\end{theorem}

As a generalization of {\rm\cite[Theorem 1.1]{shili}}, the following corollary emerges by Theorem~\ref{th1}.

\begin{cor}\ \ \label{c1} Any non-solvable $G$ has a non-nilpotent maximal subgroup of order divisible by $p$ containing the normalizer of some Sylow subgroup for any fixed prime divisor $p$ of $|G|$.
\end{cor}

Note that when $p$ is an odd prime divisor of $|G|$, arguing as in proof of Theorem~\ref{th1}, we can trivially generalize Theorem~\ref{th1} to the case of the coprime action of groups. We omit its proof here.

\begin{theorem}\ \ \label{th2} Let $A$ and $G$ be groups such that $A$ acts coprimely on $G$ by automorphisms. Assume that $p$ is any fixed odd prime divisor of $|G|$. If
every maximal $A$-invariant subgroup of $G$ of order divisible by $p$ containing the normalizer of some $A$-invariant Sylow subgroup is nilpotent. Then

$(1)$ $G$ is solvable;

$(2)$ $G$ has a Sylow tower;

$(3)$ $G$ is $q$-nilpotent or $q$-closed for each prime divisor $q$ of $|G|$.
\end{theorem}

Let $A$ and $G$ be groups such that $A$ acts coprimely on $G$ by automorphisms. Then $G$ is said to be an $A$-simple group if $G$ has no non-trivial $A$-invariant normal subgroups.
According to Theorems~\ref{th1} and ~\ref{th2}, we pose the following question.

\begin{question}\ \ \  Let $A$ and $G$ be groups such that $A$ acts coprimely on $G$ by automorphisms. If $G$ is an $A$-simple group in which a nilpotent maximal $A$-invariant subgroup is
its Sylow $2$-subgroup and every maximal $A$-invariant subgroup of even order containing the normalizer of some $A$-invariant Sylow subgroup is nilpotent,
can we ensure that $G$ is solvable?
\end{question}

\section{Necessary Lemmas}\label{s2}

\begin{lemma} {\rm\cite[Theorem 10.4.2]{rob}}\ \ \label{l1} If $G$
has a nilpotent maximal subgroup of odd order, then $G$ is
solvable.
\end{lemma}

\begin{lemma} {\rm\cite{rose}}\ \  \label{lx} If a non-abelian simple
group $G$ has a nilpotent maximal subgroup of even order, then
$G\cong{\rm PSL}_2(p)$, where $p$ is a Fermat prime or Mersenne prime
with $p\geq 17$.
\end{lemma}

\begin{lemma} {\rm\cite[Lemma 9]{wie}}\ \ \label{l2} Let $H$ be a nilpotent Hall-subgroup of a group $G$ that is not
a Sylow subgroup of $G$. If for each prime divisor $p$ of $|H|$, assume $P\in{\rm Syl}_p(H)$ we always
have $N_G(P)=H$, then there exists a normal subgroup $K$ of $G$ such that $G=KH$ and $K\cap H=1$.
\end{lemma}

\section{Proof of Theorem~\ref{th1}}\label{s3}

\demo {\bf Part I}. We first prove that $G$ is solvable.

Let $G$ be a counterexample of the smallest order. Assume that $N$ is any non-trivial normal subgroup of $G$. It is easy to see that
the hypothesis holds for the quotient group $G/N$. If $N$ is solvable. Then $G$ is solvable since $G/N$ is also solvable by the minimality of $G$,
a contradiction. Therefore, $N$ is non-solvable and $G$ has no non-trivial solvable normal subgroups.

Let $D$ be a minimal normal subgroup of $G$. Since $D$ is non-solvable, $|D|$ has at least three distinct prime divisors by $p^{\alpha}q^{\beta}$-theorem.
Let $q$ be an odd prime divisor of $|D|$. Assume $Q\in{\rm Syl}_q(G)$. Then $Q\cap D\in{\rm Syl}_q(D)$. Note that $Q\cap D\ntrianglelefteq G$, which implies
that $N_G(Q\cap D)<G$. It is obvious that $N_G(Q)\leq N_G(Q\cap D)$. It follows that $G$ has a maximal subgroup $M$ such that $N_G(Q)\leq N_G(Q\cap D)\leq M$. By Frattini's argument,
$G=N_G(Q\cap D)D=MD$.

For $M$, we divide it into two cases.

Case $(i)$: Suppose that $M$ is nilpotent. Since $G$ has no non-trivial solvable normal subgroups, $M$ is a Hall-subgroup of $G$ and the normalizer of every Sylow subgroup of $M$ equals $M$.
Moreover, since $G$ is non-solvable and $q$ is odd, $L$ is not a Sylow subgroup of $G$ by Lemma~\ref{l1}. It follows that $G$ has a normal subgroup $K$ such that $G=K\rtimes M$.
By the maximality of $M$, $K$ is a minimal normal subgroup of $G$. Since $q\mid|M|$, one has $q\nmid|K|$. Then $D$ and $K$ are two distinct minimal normal subgroups of $G$. It follows that $D\cap K=1$
and then $G\cong G/(D\cap K)\lesssim (G/D)\times(G/K)$, which implies that $G$ is solvable since both $G/D$ and $G/K$ are solvable by the minimality of $G$, a contradiction.

Case $(ii)$: Suppose that $M$ is non-nilpotent. By the hypothesis, $p\nmid|M|$. Since $G=MD$, one has $p\mid|D|$. Let $P\in{\rm Syl}_p(G)$. Then $P\cap D\in{\rm Syl}_p(D)$.
Since $P\cap D\ntrianglelefteq G$, one has $N_G(P\cap D)<G$. Then $G$ has a maximal subgroup $L$ such that $N_G(P\cap D)\leq L$. Moreover, one has $N_G(P)\leq N_G(P\cap D)\leq L$.
By Frattini's argument, $G=N_G(P\cap D)D=LD$.
Note that $p\mid|L|$. By the hypothesis, $L$ is nilpotent.

When $p$ is odd or $p=2$ but $L>P$, arguing as $M$ is Case $(i)$ we can obtain a contradiction.

Next assume $p=2$ and $L=P\in{\rm Syl}_2(G)$. Since $G=MD=LD$, one has $MD=PD$. Then $\frac{|M||D|}{|M\cap D|}=\frac{|P||D|}{|P\cap D|}$. Note that $2\nmid|M|$. It follows that
$\frac{|P|}{|P\cap D|}=1$. Then $P=P\cap D\leq D$. Furthermore, one has $G=PD=D$. It implies that $G$ is a non-abelian simple group in which the Sylow 2-subgroup $P$ is its maximal subgroup.
By Lemma~\ref{l2}, $G\cong{\rm PSL}_2(r)$, where $r\geq 17$ is a Fermat prime or a Mersenne prime. If $2\mid \frac{r-1}{2}$. By {\rm\cite{dick}}, $G$ has a maximal subgroup $H=Z_r\rtimes Z_{\frac{r-1}{2}}$.
Note that $N_G(Z_r)=H$ and $Z_r$ is a Sylow $r$-subgroup of $G$. By the hypothesis, $H$ is nilpotent, a contradiction. If $2\nmid {\frac{r-1}{2}}$. By {\rm\cite{dick}},
$G$ has a maximal subgroup $T=D_{r-1}=Z_{\frac{r-1}{2}}\rtimes Z_2$. Note that $({\frac{r-1}{2}},r(r+1))=1$ when $2\nmid {\frac{r-1}{2}}$. Let $s$ be a prime divisor of ${\frac{r-1}{2}}$ and $S$ a Sylow $s$-subgroup of $T$.
Then $S$ is also a Sylow $s$-subgroup of $G$ and $N_G(S)=T$. By the hypothesis, $T$ is nilpotent,
also a contradiction.

All above arguments imply that the counterexample of the smallest order does not exist and so $G$ is solvable.

{\bf Part II}. We show that $G$ has normal Sylow subgroups.

Let $G$ be a counterexample of the smallest order and $N$ a minimal normal subgroup of $G$. Since $G$ is solvable, $N$ is an elementary abelian subgroup of $G$ of prime-power order.
Assume $N\cong{Z_q}^m$, $m\geq 1$. By the minimality of $G$, $G/N$ has a normal Sylow subgroup $TN/N$, where $T$ is a Sylow $t$-subgroup of $G$. Then $TN\trianglelefteq G$. By the assumption,
$T\ntrianglelefteq G$. One has $N_G(T)<G$. Then $G$ has a maximal subgroup $M$ such that $N_G(T)\leq M$. By Frattini's argument, $G=N_G(T)TN=N_G(T)N=MN$.

If $p\mid|M|$. Then $M$ is nilpotent by the hypothesis. Assume $Q_M$ is a normal Sylow $q$-subgroup of $M$. It follows that $Q_MN$ is a Sylow $q$-subgroup of $G$, a contradiction.

If $p\nmid|M|$. Then $q=p$ and $N$ is a normal Sylow $p$-subgroup of $G$, also a contradiction.

Therefore, the counterexample does not exist and then $G$ has normal Sylow subgroups.

{\bf Part III}. We show that $G$ has a Sylow tower.

Assume $|G|={r_1}^{m_1}{r_2}^{m_2}\cdots{r_s}^{m_s}$, where $r_1,\,r_2,\,\cdots,\,r_s$ are distinct prime divisors of $|G|$ and $m_i\geq 1$ for $1\leq i\leq s$. Let $R_1$ be a normal Sylow $r_1$-subgroup of $G$. For the quotient group $G/R_1$, $G/R_1$ has a normal Sylow subgroup $R_2R_1/R_1$, where $R_2$ is a Sylow $r_2$-subgroup of $G$. Then we obtain a normal subgroups series: $R_1\vartriangleleft R_2R_1\vartriangleleft\cdots\vartriangleleft
R_sR_{s-1}\cdots R_2R_1$. It follows that $G$ has a Sylow tower.

{\bf Part IV}. We prove that $G$ is $q$-nilpotent or $q$-closed for each prime divisor $q$ of $|G|$.

Let $P\in{\rm Syl}_p(G)$. For $P$, we divide it into two cases.

Case (1): $P\trianglelefteq G$. Let $P=R_1$. Assume that $R_1,\,R_2,\,\cdots,\,R_s$ are all normal Sylow subgroups of $G$. Let $E=R_1\times R_2\times\cdots\times R_s=P\times R_2\times\cdots\times R_s$. Then
$E$ is a normal Hall-subgroup of $G$. Since $G$ is solvable, $E$ has a complement $F$ in $G$ by Schur-Zassenhaus's theorem. Then $G=E\rtimes F$.

Next we show that $F$ is nilpotent.

Otherwise, if $F$ is non-nilpotent. Then $F$ has a Sylow subgroup $Q$ such that $Q\ntrianglelefteq F$. It implies that $N_F(Q)<F$. Then $F$ has a maximal subgroup $F_0$ such that
$N_F(Q)\leq F_0$. Note that $G=E\rtimes F$ and $(|E|,|F|)=1$, one has $N_G(Q)=N_E(Q)N_F(Q)\leq E\rtimes F_0$, where $E\rtimes F_0$ is a maximal subgroup of $G$. Since $p\mid|E\rtimes F_0|$,
$E\rtimes F_0$ is nilpotent by the hypothesis. Then $F_0$ is nilpotent and $[E,F_0]=1$. For every $f\in F$, one has $[E,F_0]^f=[E^f,{F_0}^f]=[E,{F_0}^f]=1$. Let ${F_0}^F=\langle {F_0}^f\mid f\in F\rangle$ be the normal closure of $F_0$ in $F$. Then $[E,{F_0}^F]=1$. Since $F_0\leq {F_0}^F\leq F$ and $F_0$ is a non-normal maximal subgroup of $F$, one has ${F_0}^F=F$. It follows that $[E,F]=1$. Arguing as above $F_0$, we have that every maximal
subgroup of $F$ containing the normalizer of some Sylow subgroup is nilpotent. Then $F$ has normal Sylow subgroups by Part II. Let $B$ be a normal Sylow subgroup of $F$. Since $[E,F]=1$, one has $[E,B]=1$.
It follows that $G=E\rtimes F\leq N_G(B)$. Then $B\trianglelefteq G$, this contradicts $B\ntrianglelefteq G$ by our assumption.

Therefore, $F$ is nilpotent.

It follows that $G$ is $q$-closed when $q\mid|E|$ and $G$ is $q$-nilpotent when $q\mid|F|$.

Case (2): $P\ntrianglelefteq G$. Let $T_1,\,T_2,\,\cdots,\,T_s$ be all normal Sylow subgroups of $G$. Assume $L=T_1\times T_2\times\cdots\times T_s$. Then $L$ is a normal Hall-subgroup of $G$. Let $M$ be
a complement of $L$ in $G$. Then $G=L\rtimes M$. Note that $p\mid|M|$. Let $M_0$ be any maximal subgroup of $M$ of order divisible by $p$ containing the normalizer of some Sylow subgroup. Then
$L\rtimes M_0$ is a maximal subgroup of $G$ of order divisible by $p$ containing the normalizer of some Sylow subgroup. By the hypothesis, $L\rtimes M_0$ is nilpotent. Then $M_0$ is nilpotent. That is,
$M$ also satisfies the hypothesis of theorem. By Part II, $M$ has a normal Sylow subgroup $C$. However, $C\ntrianglelefteq G$. Then $M\leq N_G(C)<G$. Let $S$ be a maximal subgroup of $G$ such that $M\leq N_G(C)\leq S$.
Since $p\mid|S|$, $S$ is nilpotent by the hypothesis. It follows that $M$ is nilpotent.

Therefore, $G$ is $q$-closed when $q\mid|L|$ and $G$ is $q$-nilpotent when $q\mid|M|$.\hfill\qed

\end{document}